\documentclass[12pt,a4paper,twoside,leqno]{article}
\usepackage{epsfig, cite}
\usepackage{soul}
\usepackage{graphpap,latexsym,epsf}
\usepackage{graphics}
\usepackage{algorithmic}
\usepackage{subfig}
\usepackage{url}
\usepackage[colorinlistoftodos,prependcaption,textsize=tiny]{todonotes}
\usepackage[linesnumbered,ruled]{algorithm2e}
\usepackage[T1]{fontenc}
\def\betti #1{{\color{red}#1}}

\def\massi #1{{\color{blue3}#1}}

\def\ferdi #1{{\color{green1}#1}}

\def\ele#1{{\color{blue1}#1}}

\def\alex#1{{\color{magenta}#1}}

\def\diet#1{{\color{grey1}#1}}

\def\massinew #1{{\color{green2}#1}}

\let\betti\relax
\let\massi\relax
\let\ferdi\relax
\let\ele\relax
\let\alex\relax
\let\massinew\relax

\usepackage{amsthm,amsmath,amssymb}
\usepackage{titlesec}
\titleformat{\paragraph}
       {\normalfont\itshape}{\thesubsection}{}{}
\newcommand{\boxi}{\boldsymbol{\xi}}
\newcommand{\bonu}{\boldsymbol{\nu}}
\newcommand{\boeta}{\boldsymbol{\eta}}
\newcommand{\bozeta}{\boldsymbol{\zeta}}
\newcommand{\boSigma}{\boldsymbol{\Sigma}}

\definecolor{green1}{rgb}{0.1,0.65,0}
\definecolor{green2}{rgb}{0.3,0.5,0}
\definecolor{blue1}{rgb}{0.14,0.6,1.0}
\definecolor{blue3}{rgb}{0.05,0.05,0.5}
\definecolor{whit}{rgb}{1,1,1}
\definecolor{grey1}{rgb}{0.3,0.3,0.3}

\def\trait #1 #2 #3 {\vrule width #1pt height #2pt depth #3pt}
\def\fin{
      \trait .3 5 0
      \trait 5 .3 0
      \kern-5pt
      \trait 5 5 -4.7
      \trait 0.3 5 0
\medskip}
\newtheorem{teor}{Theorem}[section]
\newtheorem{defin}[teor]{Definition}
\newtheorem{lemm}[teor]{Lemma}
\newtheorem{osse}[teor]{Remark}
\newtheorem{prop}[teor]{Proposition}
\newtheorem{defi}[teor]{Definition}
\newtheorem{coro}[teor]{Corollary}
\newtheorem{prob}[teor]{Problem}
\newtheorem{hypo}[teor]{Hypothesis}

\newcommand{\bele}{\begin{lemm}\begin{sl}}
\newcommand{\enle}{\end{sl}\end{lemm}}
\newcommand{\bedef}{\begin{defi}\begin{sl}}
\newcommand{\eddef}{\end{sl}\end{defi}}
\newcommand{\bete}{\begin{teor}\begin{sl}}
\newcommand{\ente}{\end{sl}\end{teor}}
\newcommand{\beos}{\begin{osse}\begin{rm}}
\newcommand{\eddos}{\end{rm}\end{osse}}
\newcommand{\bepr}{\begin{prop}\begin{sl}}
\newcommand{\empr}{\end{sl}\end{prop}}
\newcommand{\bepro}{\begin{prob}\begin{rm}}
\newcommand{\empro}{\end{rm}\end{prob}}
\newcommand{\bede}{\begin{defin}\begin{sl}}
\newcommand{\edde}{\end{sl}\end{defin}}
\newcommand{\beco}{\begin{coro}\begin{sl}}
\newcommand{\enco}{\end{sl}\end{coro}}
\newcommand{\behy}{\begin{hypo}\begin{sl}}
\newcommand{\enhy}{\end{sl}\end{hypo}}
\newcommand{\beq}{\begin{equation}} 
\newcommand{\eeq}{\end{equation}}
\newcommand{\beqa}{\begin{eqnarray}}
\newcommand{\eeqa}{\end{eqnarray}}

\newcommand{\duav}[1]{\langle{#1}\rangle}

\newcommand{\thspace}{\hspace{3mm}}

\newcommand{\RR}{\mathbb{R}}
\newcommand{\erre}{\mathbb{R}}
\newcommand{\NN}{\mathbb{N}}
\newcommand{\Kt}{\mathbb{K}}

\newcommand{\dx}{\,{\rm d}x}

\def\qed{\ifmmode 
  \else \leavevmode\unskip\penalty9999 \hbox{}\nobreak\hfill
  \fi
  \quad\hbox{\hskip.5em\vrule width.4em height.6em depth.05em\hskip.1em}}
\def\endproofsym{\qed}
\renewenvironment{proof}[1][Proof]{\trivlist\item[\hskip\labelsep{\hskip0pt
    {\normalfont\scshape#1.}\hskip .321429\parindent}]\ignorespaces}
{\endproofsym\endtrivlist}
\def\endnobox{\def\endproofsym{}\end{proof}\def\endproofsym{\qed}}

\newcommand{\no}{\nonumber}

\newcommand{\doma}{\partial\Omega}
\newcommand{\oma}{\Omega}

\newcommand{\ioma}{\int_\Omega}

\newcommand{\calW}{\mathcal{W}}

\newcommand{\PP}{\bf {(CP)}}

\newcommand{\nf}{\bf n}
\newcommand{\bu}{{\bf u}}
\newcommand{\ub}{{\bf u}}
\newcommand{\Ub}{{\bf U}}
\newcommand{\bv}{{\bf v}}
\newcommand{\bg}{{\bf g}}
\newcommand{\bdf}{{\bf f}}

\DeclareMathOperator{\dive}{div}

\def\vp{\mathbf{\varphi}}

\newfont{\ctv}{msam10}

\def\fine{\hfill\kern4pt \vrule height4pt depth0pt width4pt }

\numberwithin{equation}{section}

\begin{document}
\author{
\ferdi{Ferdinando Auricchio}\footnote{Dipartimento di Ingegneria Civile e Architettura (DICAR), 
Universit\`{a} di Pavia, and IMATI-C.N.R.,
Via Ferrata 5, I-27100 Pavia, Italy ({\tt ferdinando.auricchio@unipv.it}).},
\ele{Elena Bonetti}\footnote{Dipartimento di Matematica ``F. Enriques', ' Universit\`{a}
di Milano, and IMATI-C.N.R., Via Saldini 50, I-20133 Milano, Italy ({\tt elena.bonetti@unimi.it}).} , 
\massi{Massimo Carraturo}\footnote{Dipartimento di Ingegneria Civile e Architettura (DICAR), 
Universit\`{a} di Pavia,
Via Ferrata 5, I-27100 Pavia, Italy and Chair of Computational in Engineering, Technical University of Munich, Munich, Germany ({\tt massimo.carraturo01@universitadipavia.it}).} , 
\\
\diet{Dietmar H\"omberg} \footnote{Weierstrass Institute for Applied
Analysis and Stochastics, Mohrenstrasse~39, D-10117 Berlin,
Germany ({\tt  hoemberg@wias-berlin.de})},
\alex{Alessandro Reali}\footnote{Dipartimento di Ingegneria Civile e Architettura (DICAR), 
Universit\`{a} di Pavia, and IMATI-C.N.R.,
Via Ferrata 5, I-27100 Pavia, Italy ({\tt alessandro.reali@unipv.it})}
and 
\betti{Elisabetta Rocca}\footnote{Dipartimento di Matematica ``F. Casorati'', 
Universit\`{a} di Pavia, and IMATI-C.N.R., 
Via Ferrata 5, I-27100 Pavia, Italy ({\tt elisabetta.rocca@unipv.it}).}
}

\title{\bf {Structural multiscale topology optimization with stress constraint for additive manufacturing}
 \footnote{\ferdi{This work was partially supported by Regione Lombardia through the project "TPro.SL - Tech Profiles for Smart Living" (No. 379384) within the Smart Living program, and through the project "MADE4LO - Metal ADditivE for LOmbardy" (No. 240963) within the POR FESR 2014-2020 program.}
\massi{MC and AR have been partially supported by Fondazione Cariplo - Regione Lombardia through the project ``Verso nuovi strumenti di simulazione super veloci ed accurati basati sull'analisi isogeometrica'', within the program RST - rafforzamento.}
The financial support of the project Fondazione Cariplo-Regione Lombardia MEGAsTAR ``Matematica
d'Eccellenza in biologia ed ingegneria come acceleratore di una nuova
strateGia per l'ATtRattivit\`a dell'ateneo pavese'' is gratefully acknowledged. The paper also benefits from the support of  the GNAMPA (Gruppo Nazionale per l'Analisi Matematica, la Probabilit\`a e le loro Applicazioni) of INdAM (Istituto Nazionale di Alta Matematica) for EB and ER. {This research was supported by the Italian Ministry of Education, University and Research (MIUR): Dipartimenti di Eccellenza Program (2018--2022) - Dept. of Mathematics ``F. Casorati'', University of Pavia. A grateful acknowledgment goes to Dr. Ing. Gianluca Alaimo for his support and precious suggestions on additive manufacturing technology.} }}


\maketitle

\vspace{-9mm}

\noindent {\bf Abstract.} 
In this paper a  phase-field approach for structural topology optimization for a 3D-printing process which includes stress constraint and potentially multiple materials or multiscales is analyzed. First order necessary optimality  conditions are rigorously derived and a numerical algorithm which implements the method is presented. A sensitivity study with respect to \betti{some parameters} is conducted for a two-dimensional cantilever beam problem. Finally, a possible workflow to obtain a 3D-printed object from the numerical solutions is described and the final structure is printed using a fused deposition modeling (FDM) 3D printer.

\vspace{.4cm}

\vspace{4mm}

\noindent
{\bf AMS (MOS) subject clas\-si\-fi\-ca\-tion:}\thspace\betti{74P05, 49M05, 74B99.}


\section{Introduction}\label{intro}

Additive manufacturing (often denoted by AM), e.g. 3D printing, is nowadays recognized as \ferdi{a very} challenging subject of research, also due to the strategic position of AM technology with respect to many applications.  This innovative technology is, at the same time, disruptive, as well as widespread and transversal. Indeed, the applications cover several  fields, like architecture, medicine, surgery, dentistry, arts. AM is deeply   changing paradigms in design and industrial production in comparison with more traditional technologies, like casting, stamping, and milling.
This kind of  technology is based on the fact that  components or complete structures are constructed through sequences of {material layers} deposition and/or curing.
 The layer by layer fashion is obtained through deposition of fused material (in the  Fused Deposition Material - FDM technology) or by melting/sintering of powders (Selecting Laser Sintering - SLS  and Selective Laser Melting - SLM technologies). Hence, hardening and solidification of the material, prior to the application of the next layer, occur mainly by thermal actions and form the bulk part. 

In recent years, in the fields of engineering and materials science,  large efforts have been \ferdi{devoted} to model AM \ferdi{processes} with particular regard to single layer behaviour, concerning  interaction between  the temperature and the stress field, heat transfer, and mechanical aspects (see e.g. \cite{hodge2014, hussein2013, turner}).  Also mesoscopic models have been developed for the layer by layer fashion, for instance applying  a  lattice Boltzmann method for the treatment of free surface flows \cite{attar2011, attar2011bis}. In a different framework, some interesting results, that have been  recently achieved for modeling  epitaxial growth,  could be put into some relation with AM modeling (see, e.g.,~\cite{fonseca}).

However, it is still necessary to introduce tools able to produce simultaneously optimization of printing materials and adjustment of prototyping processes.  Thus, in the present paper we focus on a first typical problem occurring in additive manufacturing/3D printing processes: the problem of structural optimization consisting in trying to find the best way to distribute a material in order to minimize an objective functional\massi{~\cite{bendsoe_83, bendsoe_2003, sigmund_98, sigmund_2013}}. The shape of the domain is a-priori unknown, while known quantities are the applied loads as well as regions where we want to have holes or material. Our main interest is to find regions which should be filled by material in order to optimize some properties of the sample, which is mathematically translated in the optimization of a suitable objective functional (denoted by $\mathcal{J}$ in the rest of the paper). \betti{Since we are clearly in presence of a free-boundary problem, we decide here to handle it by means of the well-known phase-field method}. 

\ferdi{
With respect to previous papers in the literature, the main novelties here are twofold}. First,  we include in the functional a constraint on the stress $\boldsymbol{{\boldsymbol \sigma}}$, which should range 
in the physically elastic domain. In previous papers either such a constraint was not included or it was imposed pointwise (cf.~\cite{bs06}), but, in that case only partial results could be obtained on the minimization problem and no  optimality conditions could be rigorously determined. 

The second and most important novelty is the derivation  of  a multiscale phase-field concept by introducing a second order parameter representing the micro-scale. This concept allows  to obtain functionally graded material (FGM) structures. 

The classical approach to shape and topology optimization is using boundary variations in order to compute shape derivatives and to decrease the functional by deforming the boundary in a shape descendant way (cf., e.g., \cite{s80, sz92}).  Another possibility, especially in order to deal with the multiscale case, is to adopt the homogenization methods (cf., e.g., \cite{a02, e04}) or the level set method which has been exploited by several authors (cf., e.g.,  \cite{b03} and references therein). The phase-field approach has been already used in structural optimization by several authors (cf., e.g., \cite{bs06, t10, tnk10, wz04}), but still few analytical results are present in the literature (cf.~\cite{bgfs14, bc03, prw12}).  

\massi{A topology optimization based on homogenization and including a two-scale (micro-macro) relationship has been introduced in~\cite{kato_2014} for non linear elastic problems. Their approach decouples the analysis of the micro structure (performed using a phase field method) from the macro structural optimization routine which uses a classical SIMP approach instead. In~\cite{clausen_2016}, it is experimentally observed that topologycally optimized infill structures (e.g., lattice structures) present an improved buckling load compared to weight with respect to bulk material structures. Regarding topology optimization for FGM, a possible approach consists in using an unpenalized SIMP method (SIM) to obtain gray scale regions which can be mapped to different lattice volumes (cf.~\cite{brackett2011, cheng2015, panesar2018}). Nevertheless, all these approaches do not allow to clearly define the boundaries of the structure which have to be reconstructed in a second step and might lead to non-optimal results.}

\ferdi{In fact, the} present work aims at obtaining a 3D-printed model by means of a multi-scale phase-field topology optimization scheme, providing at the same time a complete derivation of the first order optimality conditions.
This choice turns out to be mathematically tractable. The related analysis is indeed not very different from the one contained in \cite{bgfs14} and we utilize differentiability results already obtained therein. However,  the convex set to which the two order parameters belong is different from the simplex used in \cite{bgfs14}, where a vectorial phase field variable is introduced in order to treat the case of multi-materials. Our approach incorporates the creation of graded material structures, which could again be generalized to allow for  {multi-material} graded structures. \betti{Moreover, our objective functional contains a constraint on the stress ${\boldsymbol{\sigma}}$ which was not present in \cite{bgfs14} and which will prove to be important for the application to the AM technology, especially in the case of lightweight structures with small material volume. 

The work is organized as follows. In Section~\ref{problem} the optimization problem is described. Section~\ref{main} presents the main analytical results. In Section~\ref{sec:NumRes} we first introduce a numerical algorithm implementing the method, then we discuss the results of a sensitivity study with respect to a problem parameter, and finally we describe a simple workflow to obtain a 3D printed structure using an FDM 3D printer. Finally, in Section~\ref{conclusions} we draw the main conclusions of the problem and possible further outlooks of the presented work. \betti{Notice that other sensitivity analyses and more comparisons with the single-material cases for a simplified cost functional have been performed in the recent paper \cite{abchrr18} by the same authors.}

\section{The problem}\label{problem}
Let us consider a component located in  an open bounded and connected set $\oma\subset\RR^3$, with a smooth boundary $\doma$, and let 
$\nf$ denote the outward unit normal to $\doma$. We assume that $\doma$ is decomposed into $\Gamma_D\cup\Gamma_g$, $\Gamma_D$ with positive measure. Indeed, we will prescribe Dirichlet boundary conditions on $\Gamma_D$ and non-homogeneous Neumann type prescription (corresponding to a traction applied) on the  part $\Gamma_g$. We also introduce the following notation: we denote by $H^1_D(\Omega;\erre^d):=\{\bv\in H^1(\Omega; \erre^d)\, :\, \bv={\bf 0} \quad \hbox{on }\Gamma_D\}$ and by $H(\dive,\Omega):=\{v\in L^2(\Omega; \erre^{d\times d})\,:\, \dive v\in L^2(\Omega; \erre^{d})\}$. 

As it is known, one of the main characteristic of AM technology is the possibility to construct objects with prescribed macroscopic and microscopic structure. We aim to introduce a model to get a combined optimization of the two scales of this structure: a macroscopic scale corresponding \ferdi{either} to the presence of material \ferdi{or to the presence of no material (i.e. voids)}, and a microscopic scale corresponding   to the microscopic density of the material. To this purpose we introduce a new double phase-fields model (cf.~also \cite{xb17} for similar approaches).  We let $\Omega_0, \, \Omega_1$ be two sets of positive measure contained in $\Omega$ such that $\Omega_0\cap \Omega_1=\emptyset$ .~\ele{We aim to introduce a model to get a simultaneous optimization of the two scales of this structure: a macroscopic scale corresponding to the presence of material (or voids), and a microscopic scale corresponding to the microscopic density of the material, when it is present. To this purpose we introduce a new two-scale  phase-field model, where the two phase parameters describe the presence of the material and its density. In particular, the parameter standing for the microscopic density of the material depends on the macroscopic phase parameter through an internal constraint.} Hence, we first introduce the phase variable $\vp$ to denote a macro-scale parameter, meaning that in the regions of $\Omega$ where $\vp=1$ we have the presence of the material, while when $\vp=0$ we have voids (\ferdi{no material}). In order to describe the micro-scale effects (corresponding to possible different densities of the material) we include in the model a second phase parameter which we denote by $\chi$ related to different microscopic configurations of the object. More precisely, $\chi$ is forced to belong to the interval $[0,\vp]$\ferdi{, so that in particular} $\chi$ is forced to be 0 where we have voids, i.e.~where $\vp=0$. Note that, within the phase-field approach, for $\vp$ we assume that the interface between the two phases (material and voids) is not sharp but diffuse with  small thickness $\gamma$ (see \eqref{diffuse}). The same diffuse interface is assumed for the microscopic density denoted by $\chi$.  

Now, let us introduce the optimization problem \ferdi{and} make precise the cost functional depending on the parameters $(\vp,\chi)$. As far as the cost functional corresponding to $\vp$, \ferdi{we}  first approximate the standard perimeter term
by a multiple of the so-called Ginzburg-Landau type functional 
\begin{equation}\label{diffuse}
\int_\Omega\left(\frac{\calW(\vp)}{\gamma}+\gamma\frac{|\nabla\vp|^2}{2} \right)\dx
\end{equation}
where $\calW$ is a potential function attaining global minima at $\vp=1$ (material) and $\vp=0$ (void). A typical example of $\calW$ is the standard double well potential $\calW(\vp)=(\vp-\vp^2)^2$. 
Let us notice that the potential $\mathcal{W}$ could include a linear term of the type $\int_\Omega \vp\, dx$. This would correspond to a minimization of the volume of the material we use to create the sample and so it would be compatible both with the experiments and also with the analytical assumptions we need to prescribe on $\mathcal{W}$ (cf.~{\bf (H1)}). 
\ferdi{Furthermore, to impose a} constraint on the admissible values for the two phase variables, we introduce  in the cost functional \eqref{cost} the following term 
$$ \ioma I_C(\vp,\chi)\, dx$$ 
where $I_C(\vp, \chi)$ denotes the characteristic function of the convex set 
\begin{equation}\label{defC}
C:=\{(\vp, \chi)\, :\, \vp\in [0,1], \quad \chi\in [0,\vp]\},
\end{equation}
that is 
$$
I_C(\vp, \chi)=\begin{cases} 
0 &\hbox{if }(\vp, \chi)\in C\\
+\infty &\hbox{otherwise.}
\end{cases}
$$
\ferdi{Assuming to deal with  a linear elasticity regime problem} under the assumption of small displacements, we denote by $\bu:\Omega\to \erre^d $ the displacement vector and by ${\boldsymbol \varepsilon}(\bu):=(\nabla \bu)^{sym}$ the linearized symmetric strain tensor.

Then, let us introduce the set of admissible designs
\begin{align}\label{uad}
\mathcal{U}_{ad}:=\{(\bu, {\boldsymbol \sigma}, \vp, \chi)\in H^1_D(\Omega; \erre^d)\times L^2(\Omega; \erre^{d\times d})\times (H^1(\Omega; \erre))^2\,:\, (\vp, \chi)\in \mathcal{C}_{ad}\},
\end{align}
based on the set of admissible controls
\begin{align}
\label{csad}
\mathcal{C}_{ad}:=\left\{(\vp, \chi)\in (H^1(\Omega; \erre))^2\cap C\, :\, \vp=0 \hbox{ a.e.~on } \Omega_0, \,
\vp=1 \hbox{ a.e.~on } \Omega_1, \, \int_{\Omega}\vp\, dx = m |\Omega|\right\}
\end{align}
and the convex set $C$ being defined in \eqref{defC}. Note that we have included a volume constraint demanding that only the fraction $m\in (0,1)$ of the available volume is filled by the material. 

The goal of structural topology optimization then is to find an optimal distribution of this material fraction characterized by the macro and micro phase field parameters $(\vp, \chi)$ acting as control parameters such that the resulting structure has a maximal stiffness. Since the inverse of stiffness is flexibility or compliance, we can rephrase this in terms of the following minimization problem:

\vspace{2mm}
\noindent$\PP$ \,\,\,Minimize the cost functional
\begin{align}
\label{cost}
{\cal J}(\bu, {\boldsymbol \sigma}, \vp, \chi)=&\, \kappa_1\ioma \left(\frac{\calW(\vp)}{\gamma}+\gamma\frac{|\nabla\vp|^2}{2} \right)\dx+\kappa_2\ioma \left(I_C(\vp,\chi)+\frac{|\nabla\chi|^2}{2} \right)\dx\,\\
\nonumber
&+\kappa_3\int_\Omega\vp\, (\bdf \cdot \bu )\dx +\kappa_4\int_{\Gamma_g} \bg\cdot \bu\dx+\kappa_5\ioma F({\boldsymbol \sigma})\dx\nonumber 
\end{align}
over $(\bu, {\boldsymbol \sigma}, \vp, \chi)\in \mathcal{U}_{ad}$, and 
subject to  the stress-strain state relation 
\begin{align}
\label{eqsigma}
&-\dive {\boldsymbol \sigma} =\vp \, \bdf \quad\hbox{in }\Omega\\
\label{bousigma}
&{\boldsymbol \sigma}\cdot \nf=\bg\quad \hbox{on }\Gamma_g\\
\label{su}
&{\boldsymbol \sigma}=\Kt(\vp, \chi){\boldsymbol \varepsilon}(\bu)\quad\hbox{in }\Omega
\end{align} 
where the $\bdf\in L^2(\Omega; \erre^d)$ is a vector volume force, $\bg\in L^2(\Gamma_g;\erre^d)$ denotes a boundary traction acting on the structure, and $\Kt$ stands for the symmetric, positive definite elasticity tensor. A possible  example of $\Kt(\vp,\chi)$ is the following interpolation matrix
\begin{equation}\label{tensorel}
\Kt(\vp,\chi)=\Kt_M(\chi)\vp^3 +\Kt_V(\chi)(1-\vp)^3,
\end{equation}
where, in order to be compatible with a possible sharp interface limit (as $\gamma\to0$), we can choose 
$\Kt_V=\gamma^2\tilde\Kt_V$,
where $\tilde\Kt_V$ denotes a fixed elasticity tensor (cf.~\cite{bgfs14}) and, e.g.,~{$\Kt_M(\chi)=\tilde\Kt_V(\chi)=\Kt_A\chi+\frac{1}{\beta}\Kt_A(1-\chi)$, with $\beta\in (0,1)$. Even if FGM are intrinsecally heterogeneous, the assumption of asymptotic homogenization can be assumed within the structure (cf.~\cite{cheng2019}). Since an experimental validations of the numerical results goes beyond the scope of this work, we assume a simple linear interpolation for the material properties, but more complex material models can be directly employed within this general framework.} In this way the object's topology is defined by the parameter $\vp$, while the stiffness of the material continuously varies according to the distribution of the parameter $\chi$. \ferdi{Note that the equations \eqref{eqsigma}-\eqref{su} correspond to the quasi-static momentum balance equation combined with Dirichlet and Neumann boundary conditions.}

\beos 
\ferdi{
Let us note that we could encompass the case of a multi-material graded structure by assuming the field variable  $\vp$ to be replaced by a vector ${\boldsymbol \vp}$. In this case, we have to rewrite the relation between $\chi$ and the new ${\boldsymbol \vp}$ (depending on the physical problem we are considering) and consequently introduce in the free energy a new convex set in place of $C$ in \eqref{defC}. Actually, for the sake of simplicity, but without loss of generality, we restrict ourselves to the case of a scalar variable $\chi$.}
\eddos

\beos
Let us point out that in the cost functional \eqref{cost}, we include the last term in order to possibly account for the stress constraint, which naturally appears in applications for example in structural engineering problems \ferdi{where we want the stress not to exceed some material dependent threshold.}
{In the ideal case, we would like to impose a maximum stress ratio based on a given stress criterion  (e.g., von Mises, Tresca, Hill, ...), such that
\begin{equation}
\sigma_{max}=\text{max}\left(\dfrac{\sigma_e}{\sigma_y}\right), 
\label{eq:maxFunct}
\end{equation}
where $\sigma_e$ is the equivalent stress depending on the chosen criterion and $\sigma_y$ the material dependent yield stress. Since this function is not differentiable, a very popular solution in the literature of topology optimization with stress constraints (cf., e.g.~\cite{Le2009,Lee2012,Zhou2017}) is to employ the $p$-norm function defined as
$$ \sigma_{PN}=\left( \int_{\Omega}\left( \dfrac{\sigma_e}{\sigma_y} \right)^p \right)^{1/p}, $$
where the parameter $p$ controls the level of smoothness of the function, with $p\rightarrow\infty$ leading to the max function of Eq.~\eqref{eq:maxFunct}. 
Finally, the function $F$ can be taken as 
$$ F({\boldsymbol \sigma})=\mid\sigma_{PN}-1\mid^2. $$
}
\eddos

In the next Section~\ref{main} we state our main analytical results concerning the proof of first order optimality conditions for $\PP$. 

\section{Main results}\label{main} 

Let us first introduce some notation in order to rewrite the state system in a weak form. Given a matrix $\Kt$, we introduce the product of matrices $\mathcal{A}$, $\mathcal{B}$
$$\duav{\mathcal{A}, \mathcal{B}}_\Kt:=\int_\Omega\mathcal{A}\, : \, \Kt\mathcal{B},$$ 
where we have used the notation $\mathcal{A}\, : \,\mathcal{B}:=\sum_{i,j=1}^d \mathcal{A}_{ij}\mathcal{B}_{ij}$. 
Then, the elastic boundary value problem (\ref{eqsigma}--\ref{su}) can be rewritten in a weak formulation as 
\begin{equation}
\label{weakel}
\duav{{\boldsymbol \varepsilon(\bu)}, {\boldsymbol \varepsilon(\bv)}}_{\Kt(\vp, \chi)}=G(\bv, \vp)\quad \forall \bv\in H^1_D(\Omega; \erre^d)
\end{equation}
where $G(\bv, \vp):=\int_\Omega \vp\, \bdf\cdot\bv\, dx+\int_{\Gamma_g}\bg\cdot\bv\,$ and \ferdi{$\Kt(\vp,\chi)$ is the elasticity tensor} defined as in \eqref{tensorel}.

Now, let us consider the following assumptions on the data introduced in Section~\ref{intro}. 
\behy\label{hyp1}
Assume that  there exist  positive constants $c_w$, $\theta$, $\Theta$, $\Lambda$ such that 
\begin{itemize}
\item[$\bf (H1)$]  $\calW\in C^1(\erre)$, $\mathcal{W}\geq -c_w$
\item[$\bf (H2)$]  $\Kt_{i,j,k,l}\in C^{1,1}(\erre^2, \erre)$, $i,j,k,l\in \{1, \dots, d\}$, $\Kt_{ijkl}=\Kt_{jikl}=\Kt_{ijlk}=\Kt_{klij}$, and 
$$\theta|\mathcal{A}|^2\leq \Kt(\vp, \chi)\mathcal{A}\,:\,\mathcal{A}\leq \Theta|\mathcal{A}|^2, \quad |\partial_\vp\Kt(\vp,\chi)\mathcal{A}\,:\,\mathcal{B}|+|\partial_\chi\Kt(\vp,\chi)\mathcal{A}\,:\,\mathcal{B}|\leq \Lambda |\mathcal{A}||\mathcal{B}|\,,$$
for all symmetric matrices $\mathcal{A}$, $ \mathcal{B}\in \erre^{d\times d}\setminus \{{\bf 0}\}$ and for all $\vp\in \erre$
\item[${\bf (H3)}$] $(\bdf, \bg)\in L^2(\Omega;\erre^d)\times L^2(\Gamma_g;\erre^d)$
\item[${\bf (H4)}$] $F\in C^1(\erre^{d\times d}; \erre^+)$ is a convex function. 
\end{itemize}
\enhy
The argument we are introducing exploits the results stated in \cite{bgfs14}. \betti{A}ctually, in our case we have to deal with two state variables $(\vp,\chi)$ \betti{and with two control parameters $(\bu, {\boldsymbol{\sigma}})$}, so that the proofs have to be adapted to the vectorial case. For the sake of coherence we also use notations introduce in \betti{the} same paper. 

First, we recall a known result on the state system \eqref{weakel} (cf.~\cite[Thm.~3.1, 3.2]{bgfs14}). 
\bete\label{exiu}
For any given $(\vp,\chi)\in L^\infty(\Omega)\times L^\infty(\Omega)$, there exists a unique $(\bu, {\boldsymbol \sigma})\in H^1(\Omega; \erre^d)\times H(\dive, \Omega)$ which fulfills \eqref{weakel}  and \eqref{su}. Moreover, there exist positive constants $C_1$ and $C_2$ such that 
\begin{equation}\label{bouu}
\|(\bu, {\boldsymbol \sigma})\|_{H^1(\Omega; \erre^d)\times H(\dive, \Omega)}\leq C_1(\|\vp\|_{L^\infty(\Omega)}+\|\chi\|_{L^\infty(\Omega)}+1)
\end{equation}
and 
\begin{equation}\label{uniel}
\|\bu_1-\bu_2\|_{H^1_D(\Omega;\erre^d)}+\|{\boldsymbol \sigma}_1-{\boldsymbol \sigma}_2\|_{L^2(\Omega, \erre^{d\times d})}\leq C_2\left(\|\vp_1-\vp_2\|_{L^\infty(\Omega)}+\|\chi_1-\chi_2\|_{L^\infty(\Omega)}\right)
\end{equation}
where $C_2$ depends on the problem data and on $\|\vp_i\|_{L^\infty(\Omega)}$, $\|\chi_i\|_{L^\infty(\Omega)}$, $i=1,2$ and $(\bu_i, {\boldsymbol \sigma}_i)=\mathcal{S}(\vp_i, \chi_i)$, being $\mathcal{S}\, :(L^\infty(\Omega))^2\to H^1_D(\Omega; \erre^d)\times L^2(\Omega, \erre^{d\times d})$ defined as the solution control-to-state operator which assigns to a given control $(\vp, \chi)$ a unique state variable $(\bu, {\boldsymbol \sigma})\in H^1_D(\Omega; \erre^d)\times L^2(\Omega, \erre^{d\times d})$.  
\ente
Then, we can state our main result related to the existence of solution to Problem $\PP$ and the derivation of first order necessary optimality conditions. 
\bete\label{exi}
The problem $\PP$ has a minimizer. 
\ente
\proof 
Let us denote by ${\cal G}_{ad}:=\{(\bu, {\boldsymbol \sigma}, \vp, \chi)\in \mathcal{U}_{ad}\,:\, (\bu, {\boldsymbol \sigma}, \vp,\chi) \hbox{ fulfills } \eqref{weakel}\}$. By virtue of \eqref{weakel}  and the Hypothesis \ref{hyp1}, and taking  $\bv=\bu$ in \eqref{weakel}, we can deduce that ${\cal J}$ is bounded from below on $ {\cal G}_{ad}$, which is not empty. Thus,  the infimum of ${\cal J}$ on ${\cal G}_{ad}$ exists and we can find a minimizing sequence $\{(\bu_k, {\boldsymbol \sigma}_k, \vp_k, \chi_k)\}\subset {\cal G}_{ad}$. Moreover, using \eqref{bouu}, we obtain that 
$${\cal J}(\bu_k, {\boldsymbol \sigma}_k, \vp_k, \chi_k)\geq \delta \left(\frac{\gamma}{2}\|\nabla \vp_k\|_{L^2(\Omega)}^2+\frac{1}{2}\|\nabla \chi_k\|_{L^2(\Omega)}^2\right)-C_\delta$$
for some $\delta>0$ and $C_\delta>0$.  This inequality follows by convexity and the boundedness of $(\vp,\chi)$ (see, e.g., \eqref{defC}).  Hence, by using the fact that $\vp_k$ belong to $[0,1]$ (cf. \eqref{defC}) for all $k\in \NN$ and by means of Poincar\'e inequality 
we obtain that the sequence $\{\vp_k\}$ is bounded in $H^1(\Omega)\cap L^\infty(\Omega)$. The same can be deduced for $\chi_k$, which is uniformly bounded, too. 
Hence, by Theorem~\ref{exiu},  we have that also the sequences of $\{(\bu_k, {\boldsymbol \sigma}_k)\}$ of corresponding states are bounded in $H^1_D(\Omega;\erre^d)\times H(\dive, \Omega)$ and that there exists, by compactness, some $(\bar\bu, \bar{\boldsymbol \sigma}, \bar\vp, \bar\chi)\in H^1_D(\Omega; \erre^d)\times H(\dive, \Omega))\times (H^1(\Omega; \erre))^2$ such that (as $k\to\infty$) at least for subsequences
\begin{align}\label{cuk}
&\bu_k\to\bar\bu\quad\hbox{weakly in } H^1_D(\Omega; \erre^d)\quad\hbox{and strongly in }L^2(\Omega;\erre^d)\\
\label{csigmak}
&{\boldsymbol \sigma}_k\to \bar {\boldsymbol \sigma}\quad\hbox{weakly in } L^2(\Omega;\erre^{d\times d})\\
\label{cvpk}
&\vp_k\to\bar\vp\quad\hbox{weakly in } H^1(\Omega)\quad\hbox{and strongly in }L^2(\Omega)\\
\label{cchik}
&\chi_k\to \bar\chi \quad\hbox{weakly in } H^1(\Omega)\quad\hbox{and strongly in }L^2(\Omega)\,.
\end{align}
Moreover, since the set $\mathcal{U}_{ad}$ is convex and closed  (and so also weakly closed), we also get $(\bar\bu,  \bar{\boldsymbol \sigma}, \bar\vp, \bar\chi)\in \mathcal{U}_{ad}$. Using $({\bf H}_1)$ and the weak lower semicontinuity of $I_C$, of norms and of $F$ (cf.~${\bf (H4)}$),  we get 
$$ \mathcal{J}(\bar\bu,  \bar{\boldsymbol \sigma}, \bar\vp, \bar\chi)\leq \lim_{k\to\infty}\mathcal{J}(\bu_k, {\boldsymbol \sigma}_k, \vp_k, \chi_k).$$
Finally, due to the fact that  $(\bu_k, {\boldsymbol \sigma}_k, \vp_k)$ fulfills \eqref{weakel} we can deduce in addition that $(\bar\bu,  \bar{\boldsymbol \sigma}, \bar\vp)$ fulfills it because $\Kt(\vp_k, \chi_k){\boldsymbol \varepsilon(\bv)}$ converges strongly to $\Kt(\bar \vp, \bar \chi){\boldsymbol \varepsilon(\bv)}$ in $L^2(\Omega; \erre^{d\times d})$ and so, using \eqref{cuk}, we get 
$$ \int_\Omega \Kt(\vp_k, \chi_k){\boldsymbol \varepsilon(\bu_k)}\,:\, {\boldsymbol \varepsilon(\bv)}\, dx\to  \int_\Omega \Kt(\bar \vp, \bar\chi){\boldsymbol \varepsilon(\bar \bu)}\,:\, {\boldsymbol \varepsilon(\bv)}\, dx\,.$$
Therefore $(\bar\bu,  \bar{\boldsymbol \sigma}, \bar\vp, \bar\chi)\in \mathcal{U}_{ad}$ turns out to be  a minimizer for $\PP$. \fine

\betti{In order to deduce first order necessary optimality conditions, w}e first introduce the linearized system with respect to the variable $\phi$ and a direction $h$ in a neighborhood of $(\bar\phi,\bar\chi)$. We use the notation
$$
(\boldsymbol{\xi}^h,\boldsymbol{\eta}^h)=\partial_\phi{\cal S}(\bar\phi,\bar\chi)h,
$$
\betti{where $(\boxi^h,\boeta^h)$ satisfies:}
\begin{align}
&-\dive\boeta^h=fh\\
&\boeta^h\cdot{\bf n}=0\\
&\boeta^h={\mathbb K}_\phi(\bar\phi,\bar\chi)h{\boldsymbol \varepsilon(\bar \bu)}+{\mathbb K}(\bar\phi,\bar\chi){\boldsymbol \varepsilon(\boxi^h)}.
\end{align}
Here $\bar{\bf u}$ stand for the first component of ${\cal S}(\bar\phi,\bar\chi)$.
Analogously, we introduce the linearized system with respect to $\chi$ in a general direction $h$. Letting 
$$
(\bozeta^h,\bonu^h)=\partial_\chi{\cal S}(\bar\phi,\bar\chi)h,
$$
\betti{where $(\bozeta^h,\bonu^h)$ satisfies:}
\begin{align}
&-\dive\bonu^h=0 \quad \hbox{in }\Omega\\
&\bonu^h\cdot{\bf n}=0 \quad \hbox{on }\Gamma_g\\
&\bonu^h={\mathbb K}_\chi(\bar\phi,\bar\chi)h{\boldsymbol \varepsilon(\bar \bu)}+{\mathbb K}(\bar\phi,\bar\chi){\boldsymbol \varepsilon(\bozeta^h)} \quad \hbox{in }\Omega.
\end{align}

We can now reformulate the optimal control problem $\PP$ by means of the so-called reduced functional 
$$j(\vp, \chi):=\mathcal{J}(\mathcal{S}(\vp, \chi), \vp, \chi)$$
which is Fr\'echet differentiable in $(H^1(\Omega)\cap L^\infty (\Omega))^2$. This fact is a consequence of the Fr\'echet differentiability of $\mathcal{J}$ (cf.~\cite[Lemma 4.2]{bgfs14}), the differentiability of the control-to-state operator (cf.~\cite[Thm. 3.3]{bgfs14}) and a standard chain rule formula (cf.~\cite[Thm.~2.20]{t10}).
In particular, we have
$$
\partial_{\vp} j(\vp, \chi) h=\mathcal{J}_{\ub}(\ub, {\boldsymbol \sigma}, \vp,\chi) \boxi^h+ \mathcal{J}_{{\boldsymbol \sigma}}(\bu, {\boldsymbol \sigma}, \vp,\chi) \boeta^h+\mathcal{J}_{\vp}(\ub, {\boldsymbol \sigma}, \vp,\chi)
$$
and
$$
\partial_{\chi} j(\vp, \chi) h=\mathcal{J}_{\ub}(\ub, {\boldsymbol \sigma}, \vp,\chi) \bozeta^h+ \mathcal{J}_{{\boldsymbol \sigma}}(\bu, {\boldsymbol \sigma}, \vp,\chi) \bonu^h +\mathcal{J}_{\chi}(\ub, {\boldsymbol \sigma}, \vp,\chi).
$$
We can now restate the Problem $\PP$ in terms of minimizers of the cost \ferdi{functional}, i.e.,
\vspace{2mm}
\noindent$\PP_R:$ \,\,\,
\begin{align}
\label{costr}
\min_{(\vp, \chi)\in \mathcal{U}_{ad}} j(\vp, \chi). 
\end{align}

Then, in order to find the first order necessary optimality conditions, we introduce the so-called \ferdi{Lagrangian:}
\begin{align}\label{lagr}
\mathcal{L}(\ub, {\boldsymbol \sigma}, \vp,\chi, \Ub, \boSigma)= &\mathcal{J}(\ub, {\boldsymbol \sigma}, \vp,\chi)-\int_\Omega {\boldsymbol \sigma} \, :\, {\boldsymbol \varepsilon(\Ub)}\, dx+\int_\Omega \bdf\cdot (\vp \Ub)\, dx\\\no
&+\int_{\Gamma_g} \bg \cdot \Ub\, dx+\int_\Omega({\boldsymbol \sigma}-{\mathbb K}(\vp,\chi){\boldsymbol \varepsilon(\ub)})\boSigma\, \betti{dx}.
\end{align}
Thus, to get minimizers we consider the partial derivatives $\mathcal{L}_\ub$ and $\mathcal{L}_{\boldsymbol \sigma}$ in direction ${\bf h}$ and impose that they are \betti{equal to zero}. From these relations, it is \betti{straightforward}, also by definition of $\mathcal{J}$, to derive the so-called adjoint equations. 
In particular, we get 
 \begin{align}\label{adj1}
&\dive(\Kt^T(\bar \vp, \bar \chi)\boSigma)=\kappa_3\bar \varphi \bdf\quad \hbox{ a.e. in } \Omega\\
\label{adj2}
&\Kt^T(\bar \vp, \bar \chi)\boSigma \cdot {\bf n}=\kappa_4\bg \quad \hbox{ a.e. on }\Gamma_g\\
\label{adj3}
&\boSigma={\boldsymbol \varepsilon(\Ub)}-\kappa_5 F_{\boldsymbol\sigma}(\bar {\boldsymbol \sigma}) \quad \hbox{ a.e. in } \Omega.
\end{align}
 
 Note that since $(\bar\vp,\bar\chi)$ is a minimizer and  $\mathcal{S}(\bar\vp, \bar\chi)=(\bar\ub, \bar{\boldsymbol \sigma})\in H_D^1(\Omega; \RR^d)\times H(\dive,\Omega)$, $(\Ub,\boSigma)\in H_D^1(\Omega; \RR^d)\times H(\dive, \Omega)$ the \ferdi{corresponding} state and adjoint variables, by convexity arguments it follows that the following inequality holds 
\begin{align}
\label{gradient}
(\mathcal{L}_{(\vp, \chi)}(\bar\ub, \bar{\boldsymbol \sigma}, \bar\vp,\bar\chi, \bar\Ub, \bar\boSigma), (\vp, \chi)-(\bar\vp, \bar\chi))\geq 0.
\end{align}
By means of this process we end up with the following main result.

\bete\label{opti}
Let $(\bar\vp, \bar\chi)$ denote a minimizer of problem $\PP_R$ and $\mathcal{S}(\bar\vp, \bar\chi)=(\bar\ub, \bar{\boldsymbol \sigma})\in H_D^1(\Omega; \RR^d)\times H(\dive,\Omega)$, $(\Ub,\boSigma)\in H_D^1(\Omega; \RR^d)\times H(\dive, \Omega)$ the \ferdi{corresponding} state and adjoint variables. Then,  ($\bar\ub, \bar{\boldsymbol \sigma}, \bar\vp,\bar\chi, \bar\Ub, \bar\boSigma)$ fulfills the  optimality system in weak sense obtained coupling the state relations \eqref{adj1}-\eqref{adj3}
 and the gradient inequality arising from \eqref{gradient}:
\begin{align*}
&\kappa_1\int_\Omega \frac{\calW'(\bar\vp)}{\gamma}(\vp-\bar\vp)\, dx +\kappa_1\gamma \int_\Omega\nabla\bar\vp\nabla(\vp-\bar\vp)\, dx+\kappa_2 \int_\Omega\nabla\bar\chi\nabla(\chi-\bar\chi)\, dx\\
&+\int_\Omega \bdf\cdot (\bar\Ub+\kappa_3 \bar \bu)(\vp-\bar\vp)\, dx-\kappa_3\int_\Omega \Kt_\vp(\bar\vp, \bar\chi)\boSigma\, :\, {\boldsymbol \varepsilon(\bar \bu)}(\vp-\bar\vp)\, dx
\\
&
-\kappa_3\int_\Omega \Kt_\chi(\bar\vp, \bar\chi)\boSigma\, :\, {\boldsymbol \varepsilon(\bar \bu)}(\chi-\bar\chi)\, dx\geq 0
\end{align*}
for all  $(\vp, \chi)\in \mathcal{C}_{ad}$ and assuming $\kappa_4=\kappa_3$. 

\ente
The proof of this Theorem follows once \eqref{gradient} is satisfied by definition of $\mathcal{L}$ in \eqref{lagr} and recalling the cost functional \eqref{cost}. 

\section{Numerical results} \label{sec:NumRes}
In this section we present an application of the presented analytical results in the engineering practice. In the first part of this section we derive the discrete formulation  of the optimization problem (\textbf{CP}) {{neglecting the stress constraint, i.e., setting $\kappa_5=0$,}} successively, we discuss a sensitivity study of the resulting optimized 2D structure, and finally we present a possible procedure to obtain from numerical results a 3D printed FGM structure. The presented results are obtained using FEniCS~\cite{lm12}, an open source library to automate the solution of mathematical models based on differential equations. 
\betti{Further numerical results on the sensitivity with respect to other parameters and more comparisons with the single-material case can be found in \cite{abchrr18}.}
\subsection{Discrete problem formulation}
\label{ssec:DiscProb}
\paragraph{Allen-Cahn gradient flow}
\label{sssec:GradientTimeStep}
To obtain a discrete version of the problem (\textbf{CP}) we employ Allen-Cahn gradient flow approach, a steepest descendant pseudo-time stepping method. Given a fixed time-step increment $\tau$ the Allen-Cahn gradient flow leads to the following set of equations:
\begin{multline}
\dfrac{\gamma_{\vp}}{\tau}\int_{\Omega}(\vp_{n+1}-\vp_{n})(\vp-\vp_{n+1})\, dx +
\kappa_1\gamma\int_\Omega\nabla\vp_{n+1}\nabla(\vp-\vp_{n+1})\, dx -
\\
%
\kappa_3\int_\Omega \Kt_\vp(\vp_n, \chi_n){\boldsymbol \Sigma}\, :\, \boldsymbol \varepsilon(\bu_{n})(\vp-\vp_{n+1})\, dx
+\dfrac{\kappa_1}{\gamma}\int_{\Omega}\calW'(\vp_n)(\vp-\vp_{n+1})\, dx \geq 0,\,\forall\left(\vp, \chi\right)\in\mathcal{C}_{ad},
 \label{eq:GradientEqPhi}
\end{multline}
\begin{multline}
\dfrac{\gamma_{\chi}}{\tau}\int_{\Omega}(\chi_{n+1}-\chi_n)(\chi-\chi_{n+1})\text{d}\, dx +
\kappa_2\int_{\Omega}\nabla\chi_{n+1}\cdot\nabla(\chi-\chi_{n+1})\, dx -
\\
\kappa_3\int_\Omega \Kt_\chi(\vp_n, \chi_n){\boldsymbol\Sigma}\, :\, {\boldsymbol \varepsilon(\bu_{n})}(\chi-\chi_{n+1})\, dx \geq 0,\,\forall\left(\vp, \chi\right)\in\mathcal{C}_{ad}
 \label{eq:GradientEqChi}
\end{multline}
to be solved under the volume constraint:
\begin{equation}
\int_{\Omega}(\vp_{n+1}-m)\, dx = 0.
 \label{eq:GradientEqLambda}
\end{equation}
\paragraph{Finite element discretization}
\label{sssec:FEM}

We then discretize the physical domain $\Omega$ employing four triangular meshes $\mathcal{Q}_u$, $\mathcal{Q}_{\vp}$, $\mathcal{Q}_{\chi}$ and $\mathcal{Q}_{U}$, one for \ferdi{each variable} of the problem. 
At the nodes of each triangular element we interpolate, by means of piecewise linear basis functions, the corresponding variables $\mathbf{u}$, $\vp$, $\chi$, and $\mathbf{U}$ together with their variations $\mathbf{v}$, $v_{\vp}$, $v_{\chi}$ and $\mathbf{v_U}$, obtaining the following finite element expansions:
\begin{align*}
& \mathbf{u}\approx\mathbf{N_u}\tilde{\mathbf{u}}, & \mathbf{v}\approx\mathbf{N_u}\tilde{\mathbf{v}}, \\
& \vp\approx\mathbf{N}_{\vp}\tilde{\boldsymbol{\vp}}, & v_{\vp}\approx\mathbf{N}_{\vp}\tilde{\mathbf{v}}_{\vp}, \\
& \chi\approx\mathbf{N}_{\chi}\tilde{\boldsymbol{\chi}}, & v_{\chi}\approx\mathbf{N}_{\chi}\tilde{\mathbf{v}}_{\chi}, \\
& \mathbf{U}\approx\mathbf{N_U}\tilde{\mathbf{U}}, & v_{U}\approx\mathbf{N_U}\tilde{\mathbf{v}}_{\mathbf{U}}, \\
\end{align*}
where
$\mathbf{N_u},\mathbf{N}_{\vp},\mathbf{N}_{\chi},\mathbf{N_U}$ are the piecewise linear shape function vectors  \ferdi{or matrices} which interpolate the nodal degrees of freedoms $\tilde{\mathbf{u}},\tilde{\boldsymbol{\vp}},\tilde{\boldsymbol{\chi}},\tilde{\mathbf{U}}$ and their variations $\tilde{\mathbf{v}},\tilde{\mathbf{v}}_{\vp},\tilde{\mathbf{v}}_{\chi},\tilde{\mathbf{v}}_{\mathbf{U}}$. Finally, the Lagrange multiplier $\lambda$ used to constrain the volume is applied using a  constant scalar value on the domain $\Omega$.
\par
We can now write the discretized version of the optimal control problem (\textbf{CP}), as follows:
\begin{equation}
\dfrac{1}{\tau}
\begin{bmatrix}
\mathbf{0}  & \mathbf{0}  & \mathbf{0} & \mathbf{0}  & \mathbf{0}\\
\mathbf{0}  & \mathbf{0}  & \mathbf{0} & \mathbf{0}  & \mathbf{0}\\
\mathbf{0}  & \mathbf{0} & \mathbf{M}^{\vp\vp} & \mathbf{0}  & \mathbf{M}^{\vp\lambda} \\
\mathbf{0}  & \mathbf{0} & \mathbf{0} & \mathbf{M}^{\chi\chi} & \mathbf{0} \\
\mathbf{0}  & \mathbf{0} & \mathbf{M}^{\lambda\vp} & \mathbf{0} & \mathbf{0}
\end{bmatrix}
\begin{bmatrix}
\tilde{\mathbf{u}}  \\
\tilde{\mathbf{U}}  \\
\tilde{\boldsymbol{\vp}}  \\
\tilde{\boldsymbol{\chi}} \\
\tilde{\lambda}
\end{bmatrix}
+
\begin{bmatrix}
\mathbf{K}^{\mathbf{u}\mathbf{u}} & \mathbf{0}  & \mathbf{0} & \mathbf{0} \\
\mathbf{0}  & \mathbf{K}^{\mathbf{U}\mathbf{U}} & \mathbf{0}  & \mathbf{0} & \mathbf{0} \\
\mathbf{0}  & \mathbf{0} & \mathbf{K}^{\vp\vp} & \mathbf{0}  & \mathbf{0} \\
\mathbf{0}  & \mathbf{0} & \mathbf{0} &\mathbf{K}^{\chi\chi} & \mathbf{0} \\
\mathbf{0}  & \mathbf{0} & \mathbf{0} & \mathbf{0} & \mathbf{0}
\end{bmatrix}
\begin{bmatrix}
\tilde{\mathbf{u}}  \\
\tilde{\mathbf{U}}  \\
\tilde{\boldsymbol{\vp}}  \\
\tilde{\boldsymbol{\chi}} \\
\tilde{\lambda}
\end{bmatrix}
 =
\begin{bmatrix}
\mathbf{f} \\
\mathbf{F} + \mathbf{q}^{\sigma}\\
\mathbf{q}^{\vp} + \mathbf{q}^{s} + \mathbf{q}^{\psi} \\
\mathbf{q}^{\chi} + \mathbf{q}^{s\prime} \\
\mathbf{q}^{\lambda} 
\end{bmatrix}
\label{eq::discreteGradientEquation}
\end{equation}
with the matrix and vector terms defined as:
\begin{align*}
&\mathbf{K}^{\mathbf{u}\mathbf{u}} =  \int_{\Omega}\nabla\mathbf{N_u}^T\mathbb{K}\nabla\mathbf{N_u}\, d\Omega,
\\
&\mathbf{K}^{\mathbf{U}\mathbf{U}} =  \int_{\Omega}\nabla\mathbf{N_U}^T\mathbb{K}\nabla\mathbf{N_U}\, d\Omega,
\\
& \mathbf{M}^{\vp\vp}=\gamma_{\vp}\int_{\Omega}\mathbf{N}^T_{\vp}\mathbf{N}_{\vp}\, d\Omega,
\\
& \mathbf{K}^{\vp\vp}=\kappa_1\gamma_{\vp}\int_{\Omega}\nabla\mathbf{N}_{\vp}^T\nabla\mathbf{N}_{\vp}\, d\Omega,
\\
& \mathbf{K}^{\chi\chi}=\kappa_2\gamma_{\chi}\int_{\Omega}\nabla\mathbf{N}_{\chi}^T\nabla\mathbf{N}_{\chi}\, d\Omega,
\\
& \mathbf{M}^{\chi\chi}=\gamma_{\chi}\int_{\Omega}\mathbf{N}^T_{\chi}\mathbf{N}_{\chi}\, d\Omega,
\\
& \mathbf{M}^{\lambda\vp}=\tau\int_{\Omega}\mathbf{N}^T_{\lambda}\mathbf{N}_{\vp}d\Omega=\left(\mathbf{M}^{\vp\lambda}\right)^T,
\\
&\mathbf{f} =  \int_{\Gamma_N}\mathbf{N_u}^T\mathbf{g}\, d\Gamma,
\\
&\mathbf{F} =  \int_{\Gamma_N}\mathbf{N}^T_{\mathbf{U}}\mathbf{g}\, d\Gamma,
\\
& \mathbf{q}^{\sigma}=\kappa_5\int_{\Omega}\mathbf{N}^T_{\mathbf{U}}F_{\boldsymbol{\sigma}}(\boldsymbol{\sigma}_{n+1}), d\Omega,
\\
& \mathbf{q}^{\vp}=\dfrac{\gamma_{\vp}}{\tau}\int_{\Omega}\left(\mathbf{N}^T_{\vp}\mathbf{N}_{\vp}\right)\boldsymbol{\tilde{\vp}}_n\, d\Omega = \mathbf{M}^{\vp\vp}\boldsymbol{\tilde{\vp}}_n,
\\
& q^{\lambda}=\int_{\Omega}m\, d\Omega,
\\
& \mathbf{q}^{s}=\int_{\Omega}\mathbf{N}_{\vp}^T \Kt_\vp(\boldsymbol{\tilde\vp}_n, \boldsymbol{\tilde\chi}_n){\boldsymbol \Sigma_{n+1}}\, :\, {\boldsymbol \varepsilon(u_{n+1})}\, d\Omega,
\\
& \mathbf{q}^{\psi}=  \dfrac{\kappa_3}{\gamma_{\vp}} \int_{\Omega}\mathbf{N}_{\vp}^T \calW'(\boldsymbol{\tilde\vp}_n)\, d\Omega.
\\
& \mathbf{q}^{\chi}=\dfrac{\gamma_{\chi}}{\tau}\int_{\Omega}(\mathbf{N}^T_{\chi}\mathbf{N}_{\chi})\boldsymbol{\tilde{\chi}}_n\, d\Omega = \mathbf{M}^{\chi\chi}\boldsymbol{\tilde{\chi}}_n,
\\
\begin{split}
& \mathbf{q}^{s\prime}= \int_{\Omega}\mathbf{N}_{\vp}^T \Kt_\chi(\boldsymbol{\tilde\vp}_n, \boldsymbol{\tilde\chi}_n){\boldsymbol \Sigma_{n+1}}\, :\, {\boldsymbol \varepsilon(u_{n+1})}\, d\Omega.
\end{split}
\end{align*}
\paragraph{A graded material algorithm}
\label{sssec:Alg}
To obtain a topologically optimized structure with continuously varying material properties, we solve the problem in \betti{\eqref{eq::discreteGradientEquation}} employing a staggered iterative approach as described in Algorithm~\ref{alg:optAlgorithmMultiMaterial}. In fact, the linear system in \betti{\eqref{eq::discreteGradientEquation}} can be split into three linear systems which we solve separately: the state equation system
\begin{equation}
\mathbf{K}^{\mathbf{u}\mathbf{u}} 
\tilde{\mathbf{u}} 
 =
\mathbf{f},
\label{eq:LinearElasticitySystem}
\end{equation}
\massinew{the adjoint problem system
\begin{equation}
\mathbf{K}^{\mathbf{U}\mathbf{U}} 
\tilde{\mathbf{U}} 
 =
\mathbf{F}+\mathbf{q}^{\sigma},
\label{eq:AdjointSystem}
\end{equation}
}
and the phase-field system
\begin{equation}
\dfrac{1}{\tau}
\begin{bmatrix}
\mathbf{M}^{\vp\vp} & \mathbf{0}  & \mathbf{M}^{\vp\lambda} \\
\mathbf{0} & \mathbf{M}^{\chi\chi} & \mathbf{0} \\
\mathbf{M}^{\lambda\vp} & \mathbf{0} & \mathbf{0}
\end{bmatrix}
\begin{bmatrix}
\tilde{\boldsymbol{\vp}}  \\
\tilde{\boldsymbol{\chi}} \\
\tilde{\lambda}
\end{bmatrix}
+
\begin{bmatrix}
\mathbf{K}^{\vp\vp} & \mathbf{0}  & \mathbf{0} \\
\mathbf{0} &\mathbf{K}^{\chi\chi} & \mathbf{0} \\
\mathbf{0} & \mathbf{0} & \mathbf{0}
\end{bmatrix}
\begin{bmatrix}
\tilde{\boldsymbol{\vp}}  \\
\tilde{\boldsymbol{\chi}} \\
\tilde{\lambda}
\end{bmatrix}
 =
\begin{bmatrix}
\mathbf{q}^{\vp} + \mathbf{q}^{s} + \mathbf{q}^{\psi} \\
\mathbf{q}^{\chi} + \mathbf{q}^{s\prime} \\
q^{\lambda}
\end{bmatrix}.
\label{eq:phaseFieldMultiMaterialSystem}
\end{equation}
The \texttt{graded-material optimization} routine defined in Algorithm~\ref{alg:optAlgorithmMultiMaterial} presents an iterative procedure where we first solve the state equation system~\betti{\eqref{eq:LinearElasticitySystem}} to get the solution vector $\mathbf{\tilde{u}}_{n+1}$, \massinew{secondly we need to solve the adjoint system~\eqref{eq:AdjointSystem},} and finally we evaluate the phase-field system~\betti{\eqref{eq:phaseFieldMultiMaterialSystem}} to obtain the two phase-field solution vectors $\boldsymbol{\tilde{\vp}^\ast}_{n+1}$, $\boldsymbol{\tilde{\chi}^\ast}_{n+1}$  together with the Lagrange multiplier $\tilde{\lambda}_{n+1}$.
Every iteration ends calling the function $\texttt{rescale}$, as defined in Algorithm~\ref{alg:rescale} to impose the constraints on the phase-field variables $\vp$ and $\chi$ directly at the nodal values.~\massinew{The \texttt{graded-material optimization} routine is then repeated until either the maximum number of iteration ($max_{iter}$) is reached or the $L^2-$norm of both phase-field variable increment $\Delta_{\vp}=\parallel\boldsymbol{\vp}_{n+1}-\boldsymbol{\vp}_{n}\parallel_{L^2}$ and microscopic density variable increment $\Delta_{\chi}=\parallel\boldsymbol{\chi}_{n+1}-\boldsymbol{\chi}_{n}\parallel_{L^2}$} are below a given tolerance ($tol$).
\begin{algorithm}
\DontPrintSemicolon
\SetKwInOut{Input}{input}\SetKwInOut{Output}{output}
\caption{\texttt{graded-material optimization}}\label{alg:optAlgorithmMultiMaterial}
\Input{$\mathcal{Q}$, $\mathcal{Q}_{\vp}$, $\mathcal{Q}_{\chi}$, $\mathcal{Q}_{\lambda}$, $\boldsymbol{\vp}_0$, $\boldsymbol{\chi}_0$}
\Output{Optimal topology}
$\boldsymbol{\vp}_n\gets\boldsymbol{\vp}_0$\;
$\boldsymbol{\chi}_n\gets\boldsymbol{\chi}_0$\;
\While{($\Delta_\vp\geq tol$ or $\Delta_\chi\geq tol$) and $n\leq max_{iter}$}
{
		$\mathbf{\tilde{u}}_{n+1}\gets$ \texttt{solve}\eqref{eq:LinearElasticitySystem}\;
		$\mathbf{\tilde{U}}_{n+1}\gets$ \texttt{solve}\eqref{eq:AdjointSystem}\;
		$(\boldsymbol{\tilde{\vp}}^\ast_{n+1},\boldsymbol{\tilde{\chi}}^\ast_{n+1},
		{\tilde{\lambda}}_{n+1})\gets$ \texttt{solve}\eqref{eq:phaseFieldMultiMaterialSystem}\;
		$\boldsymbol{\tilde{\vp}}_{n+1}\gets$ \texttt{rescale} $\left(\boldsymbol{\tilde{\vp}}^\ast_{n+1},\left[0,1\right]\right)$\;
				$\boldsymbol{\tilde{\chi}}_{n+1}\gets$ \texttt{rescale} $\left(\boldsymbol{\tilde{\chi}}^\ast_{n+1},\left[0,\vp\right]\right)$\;
  \texttt{update}($\Delta_\vp$)\;
	    $\boldsymbol{\vp}_n\gets\boldsymbol{\vp}_{n+1}$ \;
		$\boldsymbol{\chi}_n\gets\boldsymbol{\chi}_{n+1}$ \;
}
\end{algorithm}
\begin{algorithm}
\DontPrintSemicolon
\SetKwInOut{Input}{input}\SetKwInOut{Output}{output}
\caption{\texttt{rescale}}\label{alg:rescale}
\Input{$\boldsymbol{\omega}$, $\left[a,b\right]$}
\Output{Constrained solution vector}
\ForAll{$\omega_i\in\boldsymbol{\omega}$}
{
\uIf{$\omega_i < a$}
{
	$\omega_i=a$ \;
}
\uElseIf{$\omega_i > b$}
{
		$\omega_i=b$ \;
}
\Else
{
do nothing \;
}
}
\end{algorithm}
\subsection{Optimization of a cantilever beam}
\label{ssec:CantBeam}
We now apply the implemented numerical method to solve a two-dimensional optimization problem. We choose to solve the cantilever beam problem depicted in Figure~\ref{CantileverBeamSetup}, where $a=200 mm$, $b=100 mm$, $\mathbf{g}=(0,-600)\left[N/mm\right]$, \betti{$\mathbf{f}=\mathbf{0}$} and with an upper bound for the volume filling rate $m$ equal to 0.8. \ferdi{We choose a material having Young modulus $E=12.5$GPa, Poisson coefficient $\nu=0.25$, and yield stress $\sigma_y=45$MPa (ABS plastic). We set the parameter $\beta=1/6$, $\gamma_{\vp}=0.01$, $\kappa_1=400$, \massinew{$\kappa_3=\kappa_4=\kappa_5=1$}, and $\tau=1E-6$. For our sensitivity study we decided to vary the penalty parameter $\kappa_2$ of the gradient  term $\int_{\Omega}\nabla\chi_{n+1}\cdot\nabla(\chi-\chi_{n+1})\, dx$ among three different values: 40, 4000, and 400000.
\par
Figure~\ref{SingleMaterialCantileverBeam} shows the result for the reference optimized structure obtained using a single-material, i.e., setting $\beta=1$, while in Figure~\ref{KappaChiStudyBeta6} we can observe the FGM structures for the three different values of  $\kappa_2$. From this figure it is evident the strong influence of $\kappa_2$ on the final distribution of the variable $\chi$. In particular, we can observe how a too high value of the penalty parameter $\kappa_2$ delivers a structure where the variable $\chi$ is not able to properly distribute (see Figure~\ref{KappaChi400000}), while on the other hand small values of $\kappa_2$ allows too strong oscillations and the algorithm do not converge anymore (see Figure~\ref{KappaChi40}). A reasonable choice for the parameter $\kappa_2$ seems to be the one reported in  Figure~\ref{KappaChi4000}, where the variable $\chi$ gradually vary from the baseline bulk material to regions where a lower stiffness is required.~\massinew{ In Figure~\ref{KappaChiStudyVM} we can observe that the maximum value of the von Mises stress is kept always below $\sigma_y$, fulfilling the prescribed stress constraint. The overall stress distribution is very similar in all three cases. The major difference lies in the higher maximum stress values concentrated at the left corners of the structures.}
\begin{figure}[h!]
\centering
    \includegraphics[width=0.66\textwidth]{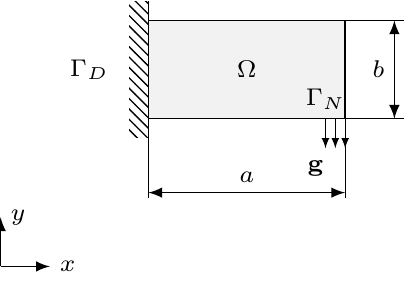}
    \caption{Cantilever beam: Problem definition.}\label{CantileverBeamSetup} 
\end{figure}
\begin{figure}[h!]
\centering
     \includegraphics[width=0.45\textwidth]{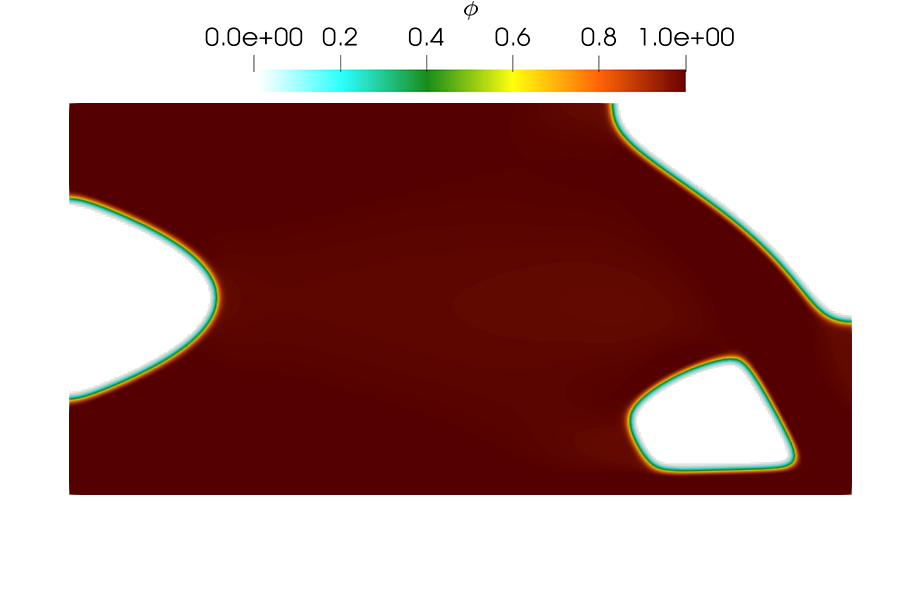}
    \caption{Cantilever beam: Reference structure obtained using a single material.}\label{SingleMaterialCantileverBeam} 
\end{figure}
\begin{figure}[h!]
\centering
  \subfloat[$\kappa_{\chi}=40$\label{KappaChi40}]
	{
     \includegraphics[width=0.3\textwidth]{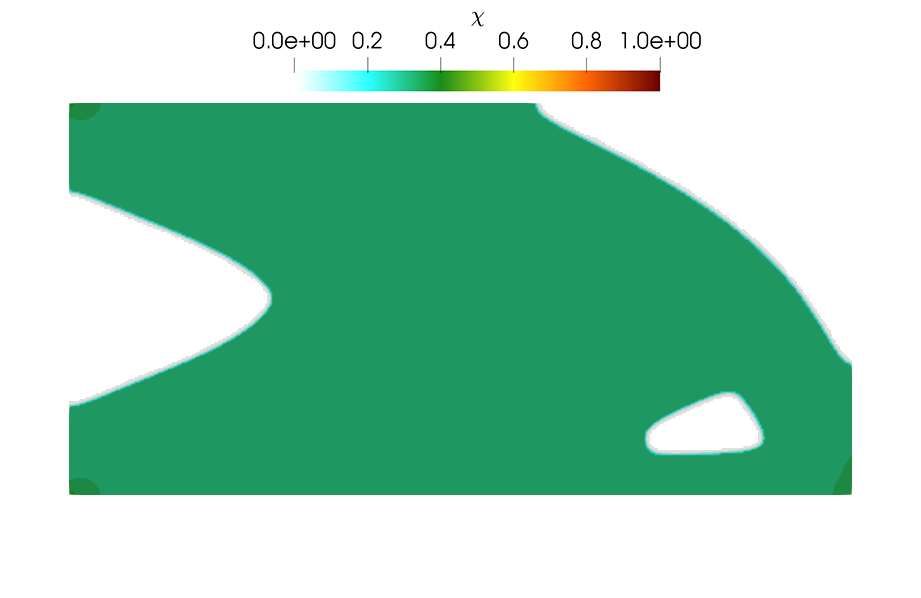}
  }
  \subfloat[$\kappa_{\chi}=4000$\label{KappaChi4000}]
	{
     \includegraphics[width=0.3\textwidth]{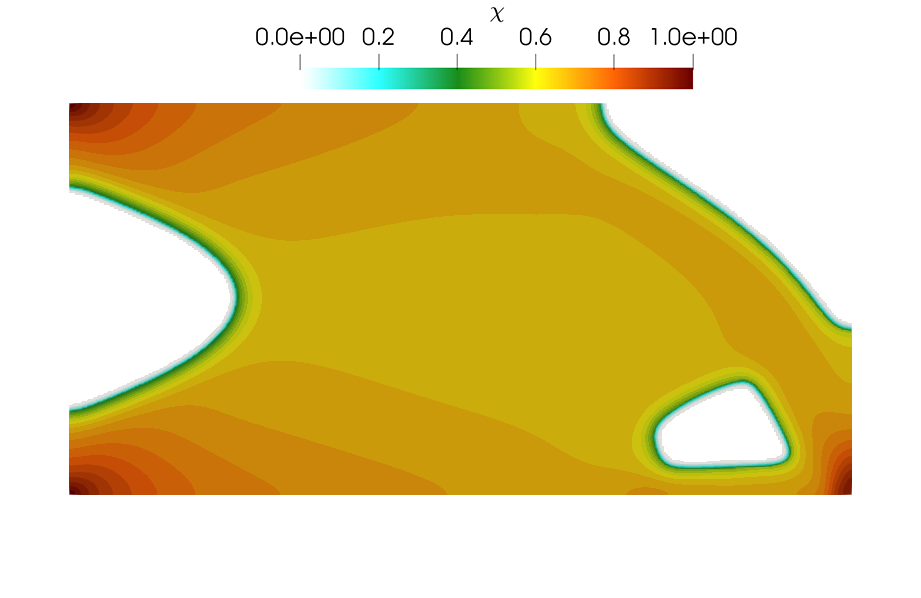}
  }
  \subfloat[$\kappa_{\chi}=400000$\label{KappaChi400000}]
	{
     \includegraphics[width=0.3\textwidth]{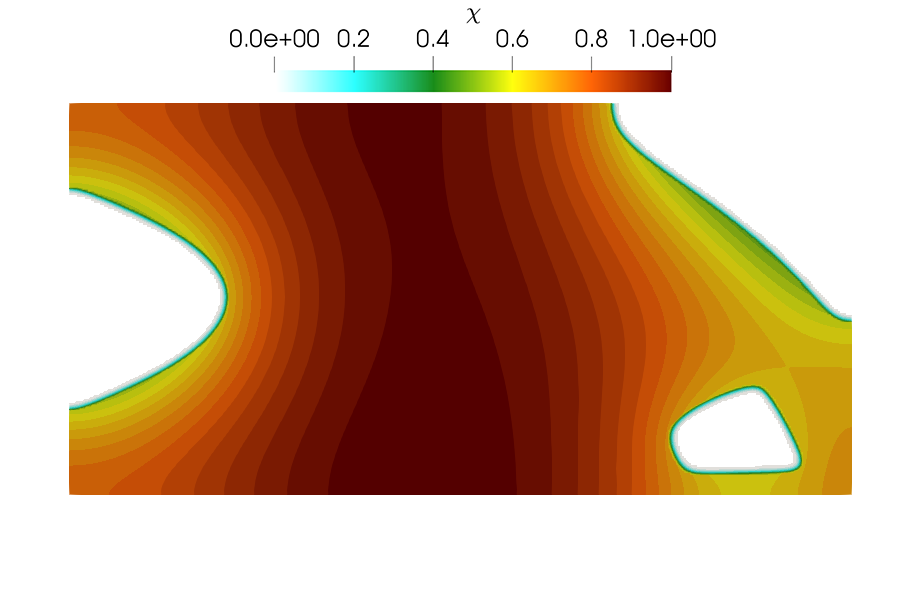}
  }
    \caption{Cantilever beam: Sensitivity study of the graded-material structure with respect to the parameter $\kappa_2$. $\chi$ value distribution.}\label{KappaChiStudyBeta6} 
\end{figure}
\begin{figure}[h!]
\centering
  \subfloat[$\kappa_{\chi}=40$\label{KappaChi40VM}]
	{
     \includegraphics[width=0.325\textwidth]{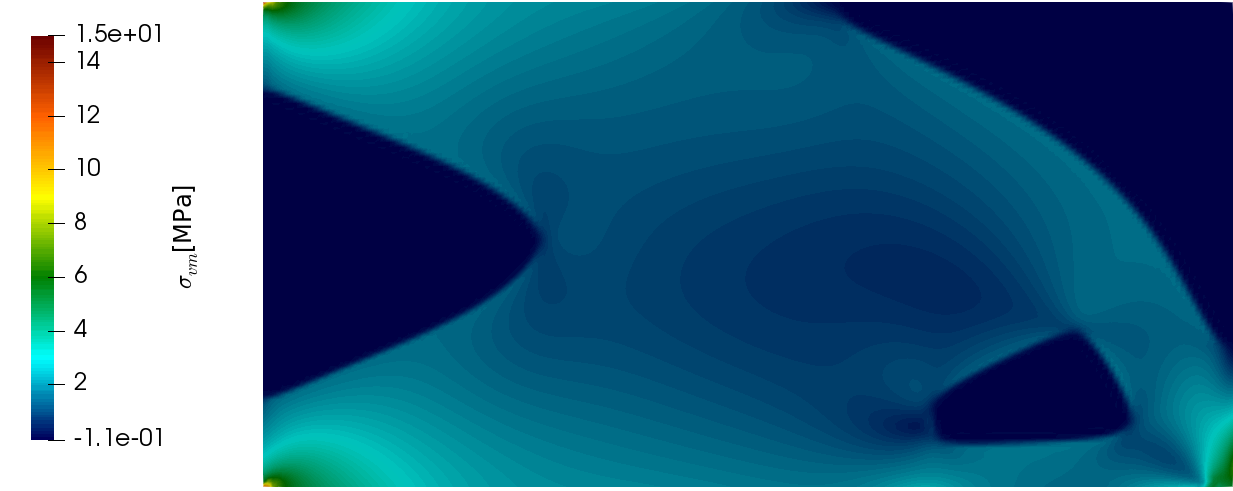}
  }
  \subfloat[$\kappa_{\chi}=4000$\label{KappaChi4000VM}]
	{
     \includegraphics[width=0.325\textwidth]{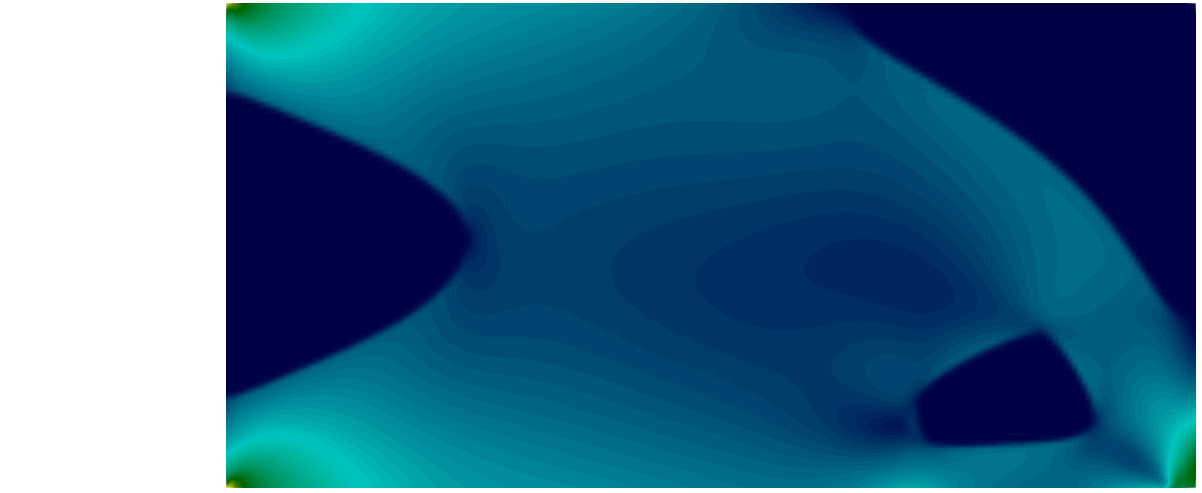}
  }
  \subfloat[$\kappa_{\chi}=400000$\label{KappaChi400000VM}]
	{
     \includegraphics[width=0.325\textwidth]{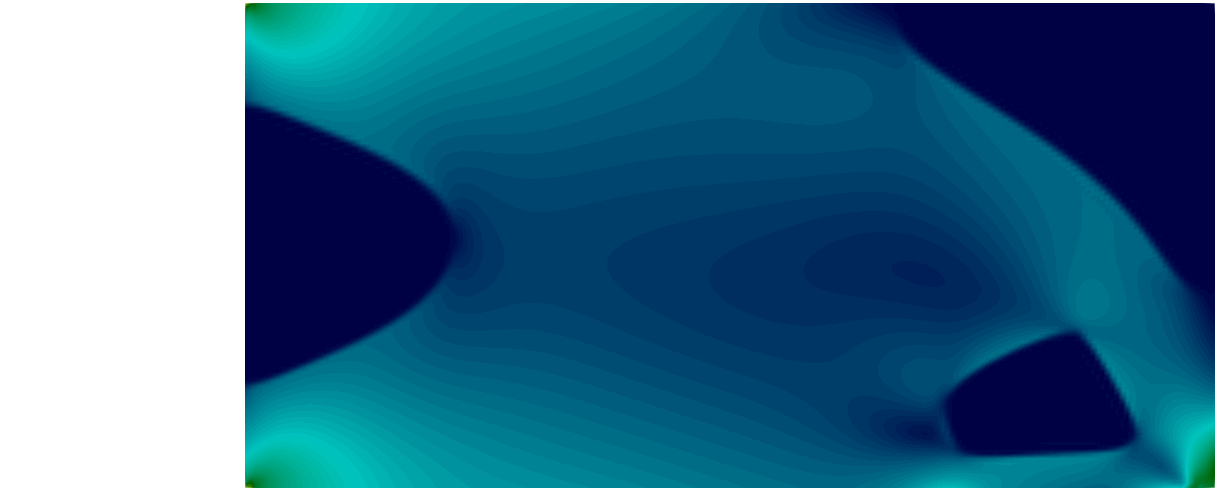}
  }
    \caption{Cantilever beam: Sensitivity study of the graded-material structure with respect to the parameter $\kappa_2$. Von Mises stress value distribution.}\label{KappaChiStudyVM} 
\end{figure}
\par
In order to estimate the total amount of material in the structure, we define a material fraction index $m_{\chi}$ as:
\begin{equation*}
m_{\chi} = \dfrac{1}{\mid\Omega\mid}\int_{\Omega}\chi d\Omega,
\end{equation*}
which can be considered as a measure of the global amount of material used to print the structure. Table~\ref{table:SensitivityGammaChi} reports the values of both the compliance and the material fraction index $m_{\chi}$ for both the single-material case and for the graded-material results. The lowest value of the compliance is achieved when a single stiffer material is used. Nevertheless, it can be observed that employing a graded-material method we are able to obtain FGM structures with a relatively low compliance using considerably less material. We want to remark here that in general the stiffness for graded-material structure do not scale linearly with the density, but it is strongly influenced by the micro-structure of the partially filled regions. This effect is not yet included with the present implementation of the method and is left to future investigations.
\begin{table}
\centering
\caption{Cantilever beam: Sensitivity study of compliance and material fraction index $m_{\chi}$ for the parameter $\kappa_2$.}
\label{table:SensitivityGammaChi}
\begin{tabular}{llll}
\hline\noalign{\smallskip}
$\kappa_2$ & compliance $\left[\dfrac{mm}{N}\right]$& $m_{\chi}$ & convergence \\
\noalign{\smallskip}\hline\noalign{\smallskip}
$40$ & $7325$ & $0.241$ & NO\\
$4000$ & $4166$ & $0.527$ & YES \\
$400000$ & $3762$ & $0.673$ & YES \\
\noalign{\smallskip}\hline
full dense material& $3130$ & $0.8$ & YES \\
\noalign{\smallskip}\hline
\end{tabular}
\end{table}

\subsection{A 3D printing workflow for topologically optimized FGM structures}
\label{ssec:3DPrint}
The topologically optimized cantilever beam of Figure~\ref{KappaChi4000} is printed using the Fused Deposition Modeling (FDM) 3D printer located at the ProtoLab~\url{http://www-4.unipv.it/3d/our-services/protolab} of the University of Pavia (see Figure~\ref{3NTR}). This machine prints a filament of thermoplastic polymer which is first heated and then extruded through a printing nozzle. The extruded filament is deployed layer by layer until the desired object is obtained (Figure~\ref{3DBeam}). Figure~\ref{PrintingPipeline} presents a very simple workflow to obtain from the numerical solution a 3D printed object. To generate a printable structure we decided to set a threshold in the $\chi$ distribution, separating the resulting structure in two regions which we print using two different plastic materials. We then extract the .STL files of the these two regions which can be now extruded and directly printed using the FDM machine.~\massinew{This extremely intuitive approach to generate printable AM structure is well suited for plastic components but it is yet not optimal. In fact, it does not allow to locally vary the material density as we observe instead in the numerical results. We are currently working on a more complex approach based on local density mapping but we leave it to forthcoming contributions.}
\begin{figure}[h!]
\centering
  \subfloat[\label{3NTR}]
	{
     \includegraphics[width=0.4\textwidth]{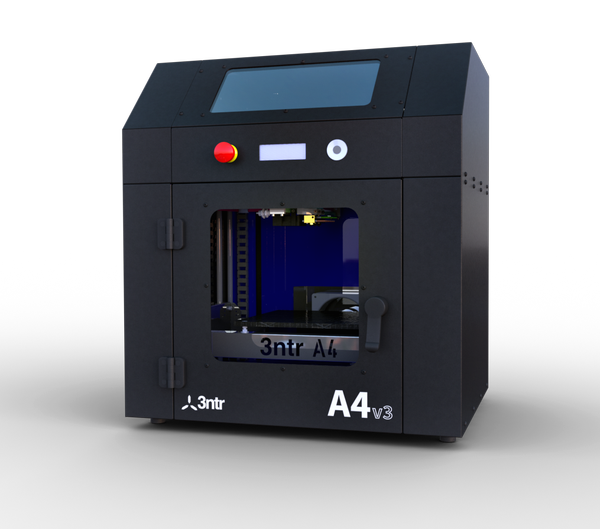}
  }
  \subfloat[\label{3DBeam}]
	{
     \includegraphics[width=0.4\textwidth]{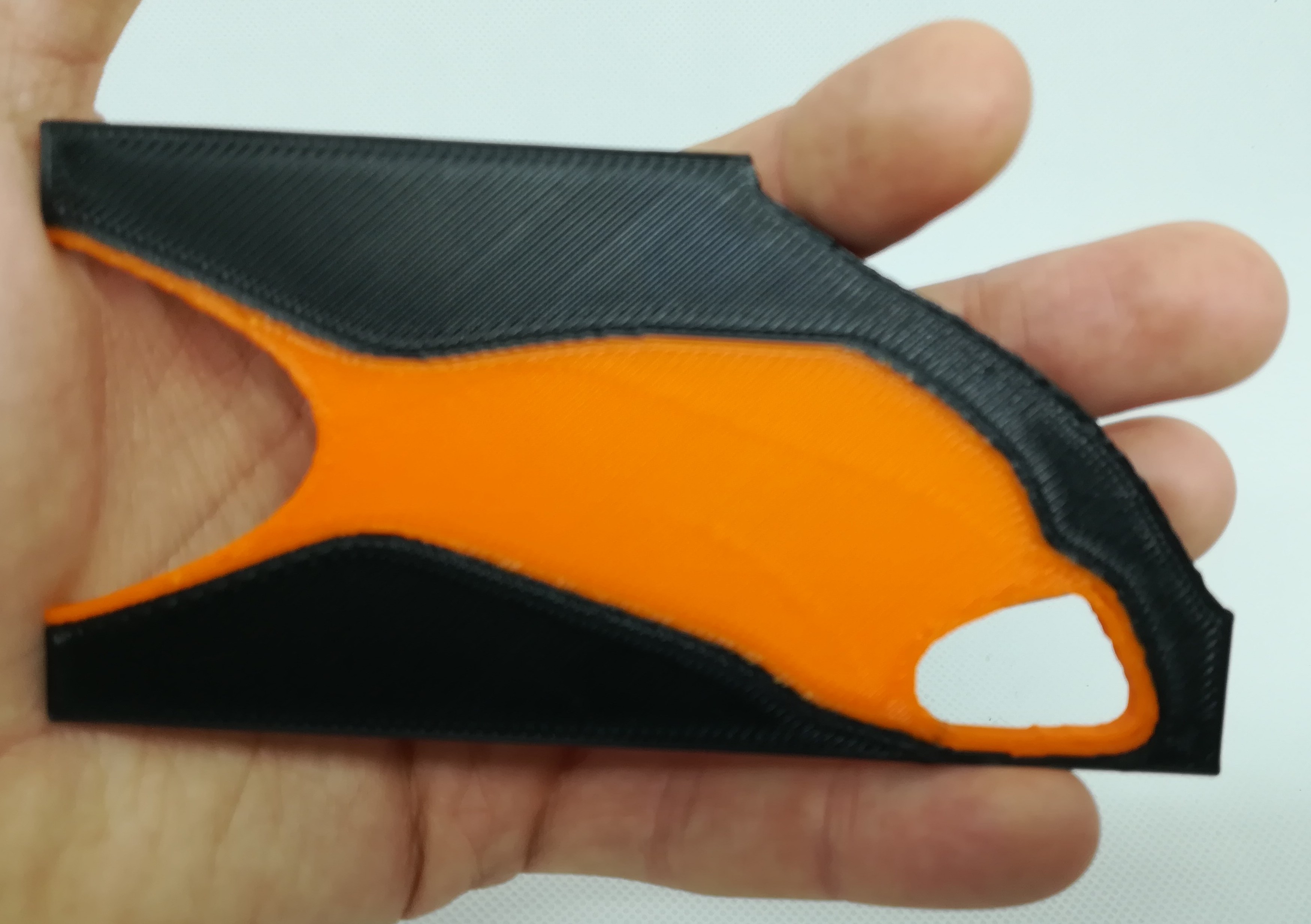}
  }
  \caption{FDM machine at ProtoLab and 3D printed cantilever beam}
\end{figure}
\begin{figure}[h!]
\centering
     \includegraphics[width=0.8\textwidth]{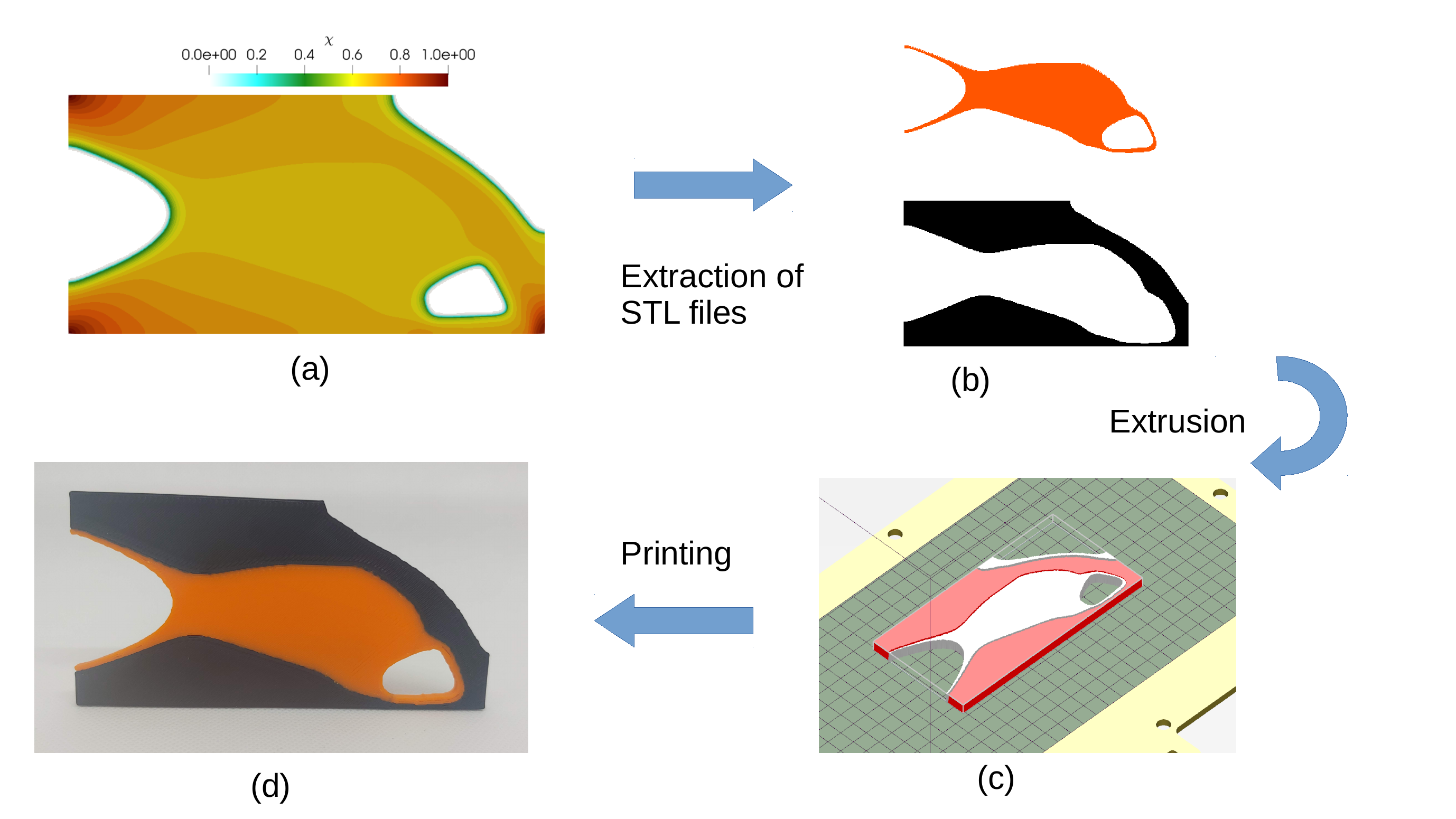}
    \caption{Description of possible workflow to obtain from a continuous $\chi$ distribution a 3D printed object: In the first step the continuous $\chi$ distribution (a) is splitted in two parts and the corresponding .STL files are generate (b), in a second step the 2D geometries are extruded to obtain a printable file (c) which can be directly sent to the FDM machine to obtain the printed structure (d).}\label{PrintingPipeline} 
\end{figure}
\section{Conclusions}\label{conclusions}
The present work analyses a phase-field approach for graded materials suitable to obtain topologically optimized structures for 3D printing processes, including stress constraints. Together with a rigorous analysis of the problem a numerical algorithm has been implemented to obtain FGM structures. A sensitivity study with respect to a problem parameter has been conducted comparing the resulting structures with a single-material reference result. Moreover, we have introduced a simple but effective workflow which from the numerical solutions leads to a 3D printed structure. Such a workflow allows us to print an optimized FGM structure using an FDM 3D printer. 
\par
As further outlooks for the present contribution we plan to investigate the influence of the microstructure on the material model and to \massinew{extend the numerical algorithm to 3D problems.

\end{document}